\documentclass[11pt,a4paper]{article}
\pdfoutput=1
\usepackage[utf8]{inputenc}
\usepackage{amsmath}
\usepackage{amsfonts}
\usepackage{amssymb}
\usepackage{amsthm}
\usepackage[a4paper,top=2.5cm,bottom=2.5cm,
            inner=3.5cm,outer=3.5cm]{geometry}
\usepackage{mathtools}
\usepackage{hyperref}
\hypersetup{colorlinks=true,allcolors=.,urlcolor=blue}
\usepackage{tikz-cd}
\usetikzlibrary{decorations.pathmorphing, arrows, decorations.markings}
\usepackage{graphicx}
\usepackage{float}
\usepackage{units} 
\usepackage{pgfplots}
\pgfplotsset{compat=1.15}
\usepackage{mathrsfs}
\usepackage[backend=biber, style=alphabetic, maxnames=5, minalphanames=3, maxalphanames=5]{biblatex}
\usepackage{ifthen}
\theoremstyle{plain}
\newtheorem{thm}{Theorem}
\numberwithin{thm}{section}
\newtheorem*{thm*}{Theorem}
\newtheorem{prp}[thm]{Proposition}
\newtheorem{lma}[thm]{Lemma}
\newtheorem{crl}[thm]{Corollary}
\newtheorem{kor}[thm]{Korollar}

\theoremstyle{definition}
\newtheorem{dfn}[thm]{Definition}
\newtheorem{exm}[thm]{Example}
\newtheorem{bsp}[thm]{Beispiel}
\newtheorem{exer}[thm]{exercise}
\newtheorem{aufg}[thm]{Aufgabe}
\newtheorem{rmk}[thm]{Remark}
\newtheorem{nta}[thm]{Notation}
\newtheorem{asp}[thm]{Assumption}

\newcommand{\theorem}[2][]{
\ifthenelse{ \equal{#1}{} }
      {\begin{thm}{#2}\end{thm}}
      {\begin{thm}[#1]{#2}\end{thm}}
}

\newcommand{\nntheorem}[2][]{
\ifthenelse{ \equal{#1}{} }
      {\begin{thm*}{#2}\end{thm*}}
      {\begin{thm*}[#1]{#2}\end{thm*}}
}

\newcommand{\proposition}[2][]{
\ifthenelse{ \equal{#1}{} }
      {\begin{prp}{#2}\end{prp}}
      {\begin{prp}[#1]{#2}\end{prp}}
}
\newcommand{\lemma}[2][]{
\ifthenelse{ \equal{#1}{} }
      {\begin{lma}{#2}\end{lma}}
      {\begin{lma}[#1]{#2}\end{lma}}
}
\newcommand{\corollary}[2][]{
\ifthenelse{ \equal{#1}{} }
      {\begin{crl}{#2}\end{crl}}
      {\begin{crl}[#1]{#2}\end{crl}}
}
\newcommand{\korollar}[2][]{
\ifthenelse{ \equal{#1}{} }
      {\begin{kor}{#2}\end{kor}}
      {\begin{kor}[#1]{#2}\end{kor}}
}

\newcommand{\definition}[2][]{
\ifthenelse{ \equal{#1}{} }
      {\begin{dfn}{#2}\end{dfn}}
      {\begin{dfn}[#1]{#2}\end{dfn}}
}
\newcommand{\example}[2][]{
\ifthenelse{ \equal{#1}{} }
      {\begin{exm}{#2}\end{exm}}
      {\begin{exm}[#1]{#2}\end{exm}}
}
\newcommand{\beispiel}[2][]{
\ifthenelse{ \equal{#1}{} }
      {\begin{bsp}{#2}\end{bsp}}
      {\begin{bsp}[#1]{#2}\end{bsp}}
}
\newcommand{\exercise}[2][]{
\ifthenelse{ \equal{#1}{} }
      {\begin{exer}{#2}\end{exer}}
      {\begin{exer}[#1]{#2}\end{exer}}
}
\newcommand{\aufgabe}[2][]{
\ifthenelse{ \equal{#1}{} }
      {\begin{aufg}{#2}\end{aufg}}
      {\begin{aufg}[#1]{#2}\end{aufg}}
}
\newcommand{\remark}[2][]{
\ifthenelse{ \equal{#1}{} }
      {\begin{rmk}{#2}\end{rmk}}
      {\begin{rmk}[#1]{#2}\end{rmk}}
}
\newcommand{\notation}[2][]{
\ifthenelse{ \equal{#1}{} }
      {\begin{nta}{#2}\end{nta}}
      {\begin{nta}[#1]{#2}\end{nta}}
}
\newcommand{\assumption}[2][]{
\ifthenelse{ \equal{#1}{} }
      {\begin{asp}{#2}\end{asp}}
      {\begin{asp}[#1]{#2}\end{asp}}
}

\renewcommand{\cal}[1]{\mathcal{#1}}
\renewcommand{\frak}[1]{\mathfrak{#1}}

\newcommand{\N}{\mathbb{N}}
\newcommand{\Z}{\mathbb{Z}}
\newcommand{\R}{\mathbb{R}}

\newcommand{\C}{\mathbb{C}}

\renewcommand{\P}{\mathbb{P}}
\newcommand{\acts}{\curvearrowright}

\newcommand{\set}[2]{\left\lbrace #1 \:\middle|\: #2 \right\rbrace}

\newcommand{\Dup}{\Delta^{\textrm{up}}}

\newcommand{\cell}{\mathrm{cell}}
\newcommand{\up}{\mathrm{up}}

\DeclareMathOperator{\tr}{tr}

\DeclareMathOperator{\Id}{Id}

\DeclareMathOperator{\Cayley}{Cayley}

\newcommand{\half}{\nicefrac{1}{2}}
\newcommand{\leftquotient}[2]{{\left.\raisebox{-.15em}{$#1$}\middle\backslash\raisebox{.05em}{$#2$}\right.}}

\bibliography{bibliography.bib} 
\author{Tim Höpfner\footnote{This work is part of the author's doctorial thesis at the University of Göttingen.}}
\title{Laplacians and Random Walks on CW Complexes}
\date{\today}

\begin{document} 
	\maketitle
	\begin{abstract} We construct random walks taking place on the $k$-cells of free $G$-CW complexes of finite type. 
These random walks define operators acting on the cellular $k$-chains that relate nicely to the (upper) cellular $k$-Laplacian. 
As an application, we use this relation to show that the Novikov-Shubin invariants of a free $G$-CW complex $X$ of finite type can be recovered from quantities related to return probabilities of the random walks on the cells of $X$. \end{abstract}
	\tableofcontents 
	\section{Introduction} 
		In this paper we generalise a classical result that relates the cellular Laplace operator acting on the $0$-chains of a locally finite CW complex $X$ to a random walk taking place on the $0$-cells of the complex.
This classical result has proven to be a useful and powerful tool, allowing us to study the Laplace operator and its spectrum by stochastic means.
Indeed, if the $1$-skeleton $X^{(1)}$ is a $d$-regular graph, and we denote the propagation operator of the uniform nearest neighbour random walk on $X^{(1)}$ by $P$, then we can interpret $P$ as an operator acting on $L^{2}C^\cell_0(X)$ and it satisfies $P = \Id - \frac{1}{d} \Delta_0$. 
Thus, we get a direct correspondence between the spectrum of $P$ and the spectrum of $\Delta_0$. 
While the assumption of pure degree makes the relation between $P$ and $\Delta_0$ particularly nice, the arguments extend without too much trouble to the case of graphs that are not of pure degree.
Recently, work of O.~Parzanchevski and R.~Rosenthal in~\cite{PaRo13} generalises this kind of results to finite simplicial complexes. 
They constructed a random walk taking place on the $k$-simplices such that the corresponding propagation operator again relates nicely to the upper $k$-Laplacian $\Dup_k = d_{k+1}^*d_{k+1}$ (see also R.~Rosenthal's paper~\cite{Ro14} and S.~Mukherjee and J.~Steenbergen's paper \cite{MuSt13} for further details). 
Using this, they deduce information about the homology of $X$ and the spectral gap of $\Delta_k$ from the random walk.
We will generalise their results further, constructing random walks that take place on free $G$-CW complexes of finite type for some group $G$ and relate the resulting random walk to $\Dup_k$. 
This builds heavily on the work of O.~Parzanchevski and R.~Rosenthal but some additional problems arise. 
In fact, we will define a family of random walks on the \textit{oriented} $k$-cells of $X$ depending on a parameter $q\in [0,1]$, such that the corresponding propagation operators $P_q$ of this random walk induces, in a natural way, an operator $B_q$ acting on the cellular $k$-chains of $X$. 
We will show that these operators satisfy the equation
$$ B_q \circ M_{1,q} = \Id - \Dup_k \circ M_{2,q}, $$
where $M_{1,q}$ and $M_{2,q}$ are multiplication operators given explicitly in terms of the local glueing information of the $k$- and $(k+1)$-cells of $X$.
Using this, we can recover spectral information of $\Dup_k$ from this random walk. 

For finite CW complexes, the classical Betti numbers measure the size of the kernels of the Laplace operators. 
Passing to infinite CW complexes, these Betti numbers tend to be infinite and provide little information. 
In these cases, we can use some extra structure, for example coming in the form of a cocompact free proper group action, to define a normalised version called the $L^2$-Betti numbers as well as the so-called Novikov-Shubin invariants that measure the content of subspaces close to the kernel of the Laplace operators.
As an application of the constructed random walk, we prove that the $k$th $L^2$-Betti number $b^{(2)}(d_{k+1}^*)$ and the $k$th Novikov-Shubin invariant $\alpha_k(X)$ (measuring the spectrum of $\Dup_k$ at zero and close to zero respectively) of a free $G$-CW complex of finite type can be described by a quantity $p_q(n)$ of the random walk that is closely related to return probabilities. 
Concretely, we show that if $X$ is regular enough so that $M_{i,q}\equiv C_{i,q}$ are constant, then for $q<1$ large enough,
$$  C_{1,q}^{-n} \left( b^{(2)}(d_{k+1}^*) + C^{-1}n^{-\alpha_k(X)/2}\right) \  \leq \  p_q(n) \ \leq\   C_{1,q}^{-n}\left(b^{(2)}(d_{k+1}^*) + Cn^{-\alpha_k(X)/2}\right). $$
For $k=0$, this recovers a classical result of N.~Th.~Varopoulos~\cite{Varop}, stating that the return probability $p$ of the uniform nearest neighbour random walk on $X^{(1)}$ satisfies $p(2n)\sim n^{-\alpha_0(X)/2}$ as $n\to \infty$.
This is a key result, as it allows us to compute the zeroth Novikov-Shubin invariant easily by using that these random walks have the same asymptotic behaviour in terms of the growth rate $N(G)$ of the group $G$, meaning that $\alpha_0(X)=N(G)$.
In contrast, there are currently only few results for Novikov-Shubin invariants in higher degrees.

Recent work in the same direction by L.~Bérnard, Y.~Chaubet, N.~V.~Dang and T.~Schick~\cite{BCDS23} uses such random walks on higher dimensional skeleta to compute $L^2$-Betti numbers and linking numbers via random walks.
	\section{Random Walks on \texorpdfstring{$k$}{k}-Cells}
		\subsection{Degree \texorpdfstring{$k$}{k}-Upper Random Walks}
			Before we define the random walk, we introduce some quantities that will be useful later on. 
These quantities capture the local structure of the CW complex described by the incidence numbers of the CW complex.
For this, let $X$ be a free $G$-CW complex of finite type\footnote{Let $G$ be a group and $X$ a CW complex with cellular left action $G\acts X$. 
Then $X$ is called a free $G$-CW complex of finite type if the projection $X\twoheadrightarrow \leftquotient{G}{X}$ is a regular covering and $\leftquotient{G}{X}$ is a finite CW complex.
In particular, the finite CW complex $\leftquotient{G}{X}$ comes with a CW structure such that there are only finitely many $k$-cells in $\leftquotient{G}{X}$, exactly one for $G$-type of $k$-cells in $X$.}. 
We denote the set of $k$-cells of $X$ by $I_k$.
On each $k$-cell $\alpha\in I_k$ we fix (arbitrarily) one of the two possible orientations.
We denote by $\alpha=\alpha_+$ the cell $\alpha$ equipped with this preferred orientation and by $-\alpha=\alpha_-$ the cell $\alpha$ with the other, reversed orientation.
Given a pair of $(k+1)$-cell $\beta\in I_{k+1}$ and $k$-cell $\alpha\in I_k$, the incidence number $[\beta:\alpha]\in \Z$ is given by the mapping degree of the map $\chi_{\beta,\alpha}$,
$$\begin{tikzcd}
\partial\beta_+ \ar[d, phantom, "\simeq" description]\ar[r, "\chi_\beta"] & X^{(k)} \ar[r, two heads] & \frac{X^{(k)}}{X^{(k)}\setminus \{\alpha\}} \ar[r, "\simeq"] & \frac{\alpha_+}{\partial\alpha_+} \ar[d, phantom, "\simeq" description]\\
S^k \ar[rrr, "\chi_{\beta,\alpha}"] & & & S^k,
\end{tikzcd}$$
where $\chi_\beta$ denotes the attaching map of $\beta$ to $X^{(k)}$.
Note that the sign of $[\beta:\alpha]$ depends on the orientations, where $[\beta_\nu:\alpha_{\nu'}] = \nu\nu'[\beta:\alpha]$ for signs $\nu,\nu'\in\{\pm\}$. 

\definition{
\label{Def_RW_ds}
Let $X$ be a free $G$-CW complex of finite type.
For $\alpha\in I_k$ and $\beta\in I_{k+1}$ we define the quantities
\begin{align}
\begin{split}
d_{+,2}(\alpha) 
&=
\sum_{\beta'\in I_{k+1}} [\beta':\alpha]^2, \\
d_+(\alpha) &= \sum_{\beta'\in I_{k+1}} |[\beta':\alpha]|, \\
d_-(\beta; \alpha) &= \begin{cases} \sum\limits_{\alpha\neq \alpha'\in I_k} |[\beta:\alpha']| &\text{ if } [\beta:\alpha]\neq 0, \\ 0 & \text{ if } [\beta:\alpha]=0, \end{cases}
\\
d_-(\alpha) &= \max_{\beta'\in I_{k+1}} d_-(\beta;\alpha),
\end{split}
\end{align}
where the maximum in the definition of $d_-(\alpha)$ exists since $d_-(\beta; \alpha)$ is invariant under the $G$-action,\footnote{$d_-(\beta; \alpha)=d_-(g.\beta; g.\alpha)$ for all $g\in G$. 
This, and the application in mind, are the reason we restrict to these complexes. The construction makes sense as long as this maximum exists.} so that it can only assume finitely many different values. Note also that these quantities are independent of the orientations chosen on $\alpha$ and $\beta$.\footnote{Since in the following, orientations on $(k+1)$-cells will never play a role, we will only take care of orientations on the $k$-cells and work with the preferred orientation on $(k+1)$-cells throughout.}
}

These quantities generalise the idea of the degree of a vertex to $k$-cells, with $d_+$ being the incoming degree and $d_-$ the (maximal) outgoing degree. 
They capture the local structure of the CW complex around the $k$-cell $\alpha\in I_k$.

We also introduce the notation
\begin{equation}
 d(\alpha, \alpha', \beta) = -[\beta:\alpha][\beta:\alpha'],
\end{equation}
measuring how well connected $\alpha$ is to $\alpha'$ along $\beta$.\footnote{
Here, we introduce an extra minus sign to mirror what happens in the case of graphs. There, random walkers can walk from a vertex $v_1$ along an (oriented) edge $e=(v_1,v_2)$ to the vertex $v_2$, where $e$ is an (outgoing) edge for $v_1$ with $[e:v_1]=-1$ and an (incoming) edge for $v_2$ with $[e:v_2]=1$, so that $[e:v_1][e:v_2]=-1$ if $e=\{v_1,v_2\}$. With this extra minus sign, the quantity $d(v_1,v_2,e) = 1$ is positive in this case.} Note that $d(\alpha,\alpha',\beta)=d(\alpha',\alpha,\beta)$. 

\remark{
Recall that on the cellular $k$-chains, the upper $k$-Laplacian $\Dup_k$ is given by the formula 
\begin{align}
\begin{split}
\label{E_Dup_Formula}
\Dup_k&\left( \sum_{\alpha\in I_k} \lambda_\alpha \cdot \alpha \right)
\\
&=\sum_{\alpha\in I_k} \left[
\sum_{\beta\in I_{k+1}} [\beta:\alpha]^2 \lambda_\alpha - \sum_{\alpha \neq \alpha'\in I_k} \sum_{\beta\in I_{k+1}} -[\beta:\alpha][\beta:\alpha'] \lambda_{\alpha'}
\right] \cdot \alpha\\
&=\sum_{\alpha\in I_k} \left[
d_{+,2}(\alpha) \lambda_\alpha - \sum_{\alpha \neq \alpha'\in I_k} \sum_{\beta\in I_{k+1}} d(\alpha,\alpha',\beta) \lambda_{\alpha'}
\right] \cdot \alpha,
\end{split}
\end{align}
where, for $k=0$, $d_{2,+}(\alpha)=\deg(\alpha)$ is the degree of the vertex $\alpha$ (ignoring loops) and $d(\alpha,\alpha',\beta)=1$ if $\beta$ is an edge from $\alpha$ to $\alpha'$ (again, ignoring loops) and zero otherwise.
Note that $\Delta^\up_k$ sees the incidence numbers between $k$- and $(k+1)$-cells in $X$.
Hence, they have to appear in the definition of our random walk.
Furthermore, both the sign and the size of the incidence numbers have to play a role.
}

\definition{
Consider the set of oriented $k$-cells
$ I_k^\pm = I_k \cup \set{\alpha_-}{\alpha \in I_k}. $
We call two oriented $k$-cells $\alpha_\nu, \alpha'_{\nu'}\in I_k^\pm$, where $\nu,\nu'\in\{\pm\}$ denote orientations, (upper) neighbours along $\beta\in I_{k+1}$, and write $\alpha_\nu \overset{\beta}{\sim} \alpha'_{\nu'}$ if 
$$\alpha_\nu\neq -\alpha'_{\nu'} \quad \text{ and } \quad d(\alpha_\nu,\alpha'_{\nu'},\beta) = \nu\nu' d(\alpha,\alpha',\beta) > 0. $$
 Note that this condition is independent of the orientation chosen on $\beta\in I_{k+1}$.
}

\definition{
The random walk $\frak R^k = \frak R^k(X)$ is given by the state space $I_k^* = I_k^\pm \cup \{\Theta\}$, where $\Theta$ is an auxiliary, absorbing state together with the following moving probabilities:
\begin{itemize}
\item The moving probabilities starting from the absorbing state $\Theta$ are given by
\begin{align*}
\P(\Theta\to \Theta) &= 1 \quad \text{and} \quad
\P(\Theta\to \alpha_\pm) = 0 \quad \text{for all } \alpha_\pm\in I_k^\pm.
\end{align*}
\item  To define the moving probabilities starting from $\alpha_\nu\in I_k^\pm$, we define first for $\alpha'_{\nu'}\in I^\pm_k$ and $\beta\in I_{k+1}$ the quantities
\begin{align*}
 \P(\alpha_\nu\nearrow \beta) &= \frac{|[\beta:\alpha_\nu]|}{d_+(\alpha)} \qquad \text{and} \qquad
 \P_\alpha(\beta\searrow \alpha'_{\nu'})= \frac{|[\beta:\alpha'_{\nu'}]|}{d_-(\alpha)}
\end{align*}
These probabilities can be seen as an intermediate step of moving first from $\alpha$ to $\beta$ and then from $\beta$ to $\alpha'$ (keeping in mind that we started at $\alpha$). 
In this sense, we define 
\begin{align*}
&\P(\alpha_\nu\xrightarrow{\beta}\alpha'_{\nu'})  \\
&= \begin{cases} \P(\alpha_\nu \nearrow \beta) \P_\alpha(\beta\searrow\alpha'_{\nu'}) 
= \frac{-[\beta:\alpha_\nu][\beta:\alpha'_{\nu'}]}{d_+(\alpha)d_-(\alpha)} 
= \frac{d(\alpha_\nu,\alpha'_{\nu'},\beta)}{d_+(\alpha)d_-(\alpha)} > 0
& \text{if } \alpha_\nu \overset{\beta}{\sim} \alpha'_{\nu'}, \\ 0 & \text{else,}
\end{cases} 
\end{align*}
the probability of moving from $\alpha$ along $\beta$ to $\alpha'$.
Recall here that $\alpha\not\sim \pm \alpha$ by definition, so that $\P(\alpha\xrightarrow{\beta} \pm \alpha) = 0$.
Finally, we set
\begin{align*}
 \P(\alpha_\nu\to\alpha'_{\nu'})  
 &= \sum_{\beta\in I_{k+1}} \P(\alpha_\nu\xrightarrow{\beta}\alpha'_{\nu'}) 
= \sum_{\stackrel{\beta\in I_{k+1}}{\alpha \overset{\beta}{\sim} \alpha'}} \frac{d(\alpha_\nu,\alpha'_{\nu'},\beta)}{d_+(\alpha)d_-(\alpha)}
\\
&= \frac{1}{{d_+(\alpha)d_-(\alpha)}}\sum_{\stackrel{\beta\in I_{k+1}}{\alpha_\nu \overset{\beta}{\sim} \alpha'_{\nu'}}} {d(\alpha_\nu,\alpha'_{\nu'},\beta)}
\end{align*}
The moving probabilities would add to one if we use $d_-(\beta;\alpha)$ in place of $d_-(\alpha)$, however it will be important later on that we can pull the factor $d_-(\alpha)^{-1}$ out of the sum as it depends only on the starting cell $\alpha$.
Consequently, their sum needs not to be one but can be smaller since we might make the denumerator bigger for some of the summands.
\item 
The complementary probability will be the probability of moving to $\Theta$, that is 
$$ \P(\alpha_\nu \to \Theta) = 1 - \sum_{\alpha'_{\nu'}\in I_k^\pm} \P(\alpha_\nu \to \alpha'_{\nu'}). $$
\end{itemize}
This defines a random walk on $I_k^\ast$.  
The propagation operator $P$ with entries 
$ P_{s,s'} = \P( s'\to s)$ for $s,s'\in I_k^*$ acts on $\ell^2(I_k^*) = \set{\sum_{s\in I_k^*}\lambda_s \cdot s}{\sum_{s\in I_k^*}|\lambda_s|^2<\infty}$ by
$$ P\left(\sum_{s\in I_k^*}\lambda_s \cdot s\right) = \sum_{s\in I_k^*}\left[\sum_{s'\in I_k^*} P_{s,s'} \lambda_{s'}\right] \cdot s = \sum_{s\in I_k^*}\left[\sum_{s'\in I_k^*} \P(s'\to s) \lambda_{s'}\right] \cdot s. $$ }

\example{
Let us consider a small example to see how we can find the moving probabilities.
For this, let $X$ be a CW complex and $\alpha\in I_k$ some $k$-cell. 
In order to find the $k$-cells we can move to from $\alpha=+\alpha$, we proceed as follows:
\begin{enumerate}
\item First we consider all $(k+1)$-cells $\beta\in I_{k+1}$ and find those, that have non-zero incidence number $[\beta:\alpha]\neq 0$. 
For this example, let us say there are two such $(k+1)$-cells $\beta_1$ and $\beta_2$ with $[\beta_1: \alpha] = 1$ and $[\beta_2:\alpha] = -2$. 
\item Then we consider all other $k$-cells $\alpha\neq \alpha'\in I_k$ that have non-zero incidence numbers with at least one of the $(k+1)$-cells above, that is $[\beta_1:\alpha']\neq 0$ or $[\beta_2:\alpha']\neq 0$. For this example, let us say there are three such $k$-cells, $\alpha_1$ with $[\beta_1:\alpha_1] = 1$, $\alpha_2$ with $[\beta_1:\alpha_2] = 2$ and $[\beta_2:\alpha_2] = 4$ and $\alpha_3$ with $[\beta_2:\alpha_3] = -2$. 
\end{enumerate}
We visualise this with the following diagram:
$$
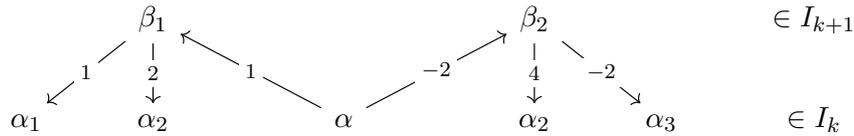
\begin{tikzcd}[ampersand replacement = \&]
\& \beta_1 \ar[dl, "1" description] \ar[d, "2" description] \& \& \& \& \beta_2 \ar[d, "4" description] \ar[dr, "-2" description] \&
 \& \in I_{k+1} \\
\alpha_1 \& \alpha_2 \& \& \alpha \ar[ull, "1" description] \ar[urr, "-2" description]  \& \& \alpha_2 \& \alpha_3
 \& \in I_{k} 
\end{tikzcd}$$
Next, we first change the orientations on the $(k+1)$-cells such that the incidence numbers with $\alpha$ are negative, so here we change the orientation on $\beta_1$.\footnote{This orientation of $\beta_1$ has no impact on the formulas in the end and is thus suppressed in the formal definition.} 
Then we change the orientations on the $\alpha_i$, $i\in\{1,2,3\}$, for each of the $\beta_j$ independently  such that the incidence numbers with the $(k+1)$-cells become positive.
This changes our diagram as follows:
$$
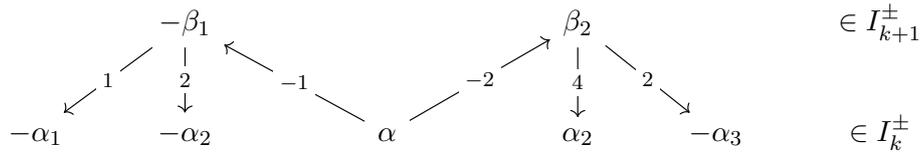
\begin{tikzcd}[ampersand replacement = \&]
\& -\beta_1 \ar[dl, "1" description] \ar[d, "2" description] \& \& \& \& \beta_2 \ar[d, "4" description] \ar[dr, "2" description] \&
 \& \in I_{k+1}^\pm \\
-\alpha_1 \& -\alpha_2 \& \& \alpha \ar[ull, "-1" description] \ar[urr, "-2" description]  \& \& \alpha_2 \& -\alpha_3 
 \& \in I_{k}^\pm
\end{tikzcd}$$
We now introduce the auxiliary state $\Theta$. 
For each of the $(k+1)$-cells, we sum the outgoing incidence numbers.
Here, we get $d_-(\beta_1;\alpha) = 1+2 = 3$ for $\beta_1$ and $d_-(\beta_2;\alpha) =  2+4 = 6$ for $\beta_2$. 
The maximum is therefore $6$, and we add connections from each of the $(k+1)$-cells to the new state $\Theta$ until the sum of outgoing edges is equal to this maximum, i.e., $-\beta_1\xrightarrow{3}\Theta$ in this case:
$$
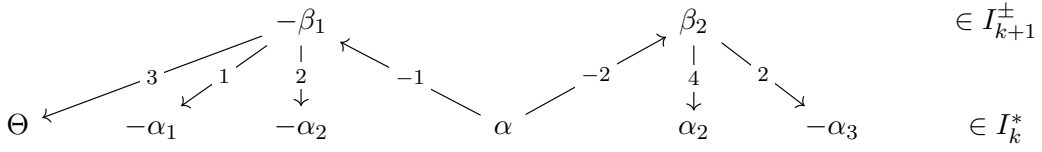
\begin{tikzcd}[ampersand replacement = \&]
\& \& -\beta_1 \ar[dll, "3" description] \ar[dl, "1" description] \ar[d, "2" description] \& \& \& \& \beta_2 \ar[d, "4" description] \ar[dr, "2" description]  \&
 \& \in I_{k+1}^\pm\\
\Theta \& -\alpha_1 \& -\alpha_2 \& \& \alpha \ar[ull, "-1" description] \ar[urr, "-2" description]  \& \& \alpha_2 \& -\alpha_3  
 \& \in I_{k}^*
\end{tikzcd}$$
The moving probabilities can now be read to be proportional to the annotations of the arrows. 
For the first intermediate step, $\P(\alpha\nearrow -\beta_1) = \nicefrac{1}{3}$ and $\P(\alpha\nearrow \beta_2) = \nicefrac{2}{3}$. 
For the second step, 
$$\P(-\beta_1\searrow \Theta) = \nicefrac{3}{6}, \qquad \P(-\beta_1\searrow -\alpha_1) = \nicefrac{1}{6} \quad \text{and} \quad \P(-\beta_1 \searrow -\alpha_2) = \nicefrac{2}{6}$$ 
and for $\beta_2$ we have 
$$\P(\beta_2\searrow \alpha_2) = \nicefrac{4}{6} \quad \text{and} \quad \P(\beta_2\searrow -\alpha_3) = \nicefrac{2}{6}.$$
The introduction of $\Theta$ guaranties that the denominator of these (unreduced) fractions is the same everywhere and the moving probabilities are proportional to the incidence numbers even if the first intermediate step leads to different $\beta_i$s.
Multiplying these accordingly, we find
\begin{alignat*}{3}
&\P(\alpha\xrightarrow{-\beta_1} -\alpha_1) = \nicefrac{1}{3}\cdot\nicefrac{1}{6}, \quad 
&& \P(\alpha\xrightarrow{-\beta_1} -\alpha_2) = \nicefrac{1}{3}\cdot\nicefrac{2}{6}, \quad
&& \P(\alpha\xrightarrow{-\beta_1} \Theta) = \nicefrac{1}{3}\cdot \nicefrac{3}{6}, \\
&\P(\alpha\xrightarrow{\beta_2} \alpha_2) = \nicefrac{2}{3}\cdot\nicefrac{4}{6}, && \P(\alpha\xrightarrow{-\beta_1} -\alpha_3) = \nicefrac{2}{3}\cdot\nicefrac{2}{6}. &&
\end{alignat*} 
Here, every oriented $k$-cell can be reached only via one $(k+1)$-cell, otherwise we would have to sum over all $(k+1)$-cells, that is $\P(\alpha\to s) = \sum_{\pm\beta\in I_{k+1}^\pm} \P(\alpha\xrightarrow{\pm\beta} s)$ for $s\in I_k^*$.
Therefore, in this example, a random walker starting at $\alpha$ has the following possible moves, with annotations denoting the probabilities:
$$
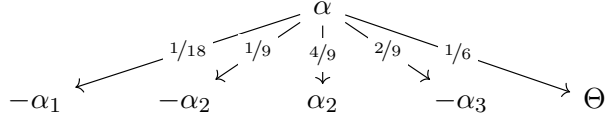
\begin{tikzcd}[ampersand replacement = \&]
\& \& \alpha \ar[dll, "\nicefrac{1}{18}" description] \ar[dl, "\nicefrac{1}{9}" description] \ar[d, "\nicefrac{4}{9}" description] \ar[dr, "\nicefrac{2}{9}" description] \ar[drr, "\nicefrac{1}{6}" description] \& \& \\
-\alpha_1 \& -\alpha_2 \& \alpha_2 \& -\alpha_3 \& \Theta
\end{tikzcd}$$
Note that the cell $\alpha_2$ can be reached with both possible orientations.
If we start at the cell $-\alpha$ with reversed orientation, we obtain the same moving probabilities as for $\alpha$, but now leading to the same cells but with flipped orientations instead. 
}

We now define an operator $B$ acting directly on the unoriented $k$-skeleton $I_k$, which is closely related to this random walk.

\definition{
We define the projection operator $T\colon \ell^2(I_k^*) \to \ell^2(I_k)$ by
$$ T\left(\sum_{s\in I_k^*}\lambda_s \cdot s\right) = \sum_{\alpha\in I_k} (\lambda_{\alpha_+} - \lambda_{\alpha_-}) \cdot \alpha $$
and the inclusion operator $I\colon \ell^2(I_k)\to \ell^2(I_k^*)$, using that $I_k=I_k^+\subset I_k^*$, by
$$ I\left( \sum_{\alpha\in I_k} \lambda_\alpha \cdot \alpha \right) = \sum_{\alpha_+\in I_k^+} \lambda_{\alpha_+}\cdot \alpha_+. $$
Lastly, we define the operator $B\colon \ell^2(I_k)\to \ell^2(I_k)$ by
$$ B\left(\sum_{\alpha\in I_k} \lambda_\alpha \cdot \alpha \right) 
= \sum_{\alpha\in I_k} \left[\sum_{\alpha\neq \alpha'\in I_k} \frac{1}{{d_+(\alpha')d_-(\alpha')}} \sum_{\beta\in I_{k+1}}  d(\alpha,\alpha',\beta) \lambda_{\alpha'}\right]\cdot \alpha.$$
With $B_{\alpha,\alpha'} = \frac{1}{{d_+(\alpha')d_-(\alpha')}} \sum_{\beta\in I_{k+1}} d(\alpha_\nu,\alpha'_{\nu'},\beta)$ for $\alpha\neq\alpha'$ and $B_{\alpha,\alpha}=0$, 
$$ B\left(\sum_{\alpha\in I_k} \lambda_\alpha \cdot \alpha \right) = \sum_{\alpha\in I_k} \left[\sum_{\alpha'\in I_k} B_{\alpha,\alpha'} \lambda_{\alpha'}\right]\cdot \alpha.$$
}
This is captured by the diagram
$$\begin{tikzcd}[row sep = 30, column sep = 30]
\ell^2(I_k^*) \ar[r, "P"]\ar[d, "T"] & \ell^2(I_k^*)\ar[d, "T"] \\
\ell^2(I_k) \ar[r, "B"]  \ar[u, bend left, dashed, "I"] & \ell^2(I_k). \ar[u, bend left, dashed, "I"] 
\end{tikzcd}$$
This operator $B$ does not describe a random walk since $B_{\alpha,\alpha'}$, the ``probability of moving from $\alpha'$ to $\alpha$'', may even be negative.
However, this operator is closely related to the random walk described by $P$.
Indeed, using the operators $T$ and $I$ we can see that $B$ describes the process that arises from the random walk if we consider a random walker arriving at a cell $\alpha_-$ (equipped with the reversed orientation) as the inverse of a random walker at $\alpha_+$ --- that is we allow random walkers at $\alpha_+$ and $\alpha_-$ to cancel each other out.\footnote{While it is not clear which, if any, physical process this operator $B$ describes, cancellation between different objects does happen in physics. 
For example, studying fermions via the Dirac equation suggests that for every particle there is a corresponding anti-particle, such as electrons and positrons. 
If they meet, they will annihilate each other.}

\lemma{
The operators $T$, $I$ and $B$ defined above satisfy the equations
$$BT = TP, \qquad B=TPI \quad \text{and} \quad B^n = TP^n I.$$ 
}
\begin{proof}
We check these equalities by direct computation. 
For $BT$ we obtain
\begin{align*}
BT\left( \sum_{s\in I_k^*} \lambda_{s}\cdot s \right)
&= B\left(\sum_{\alpha\in I_k} (\lambda_{\alpha_+} - \lambda_{\alpha_-}) \cdot \alpha\right) 
\\
&= 
\sum_{\alpha\in I_k} \left[  \sum_{\alpha'\in I_k} \frac{1}{{d_+(\alpha')d_-(\alpha')}} \sum_{\beta\in I_{k+1}} d(\alpha,\alpha',\beta) \cdot (\lambda_{\alpha'_+}-\lambda_{\alpha'_-})\right] \cdot \alpha
\end{align*}
and for $TP$ we compute (omitting the coefficient of $\Theta$ as it disappears in the next step) that 
\begin{align*}
TP&\left(\sum_{s\in I_k^*} \lambda_{s} s \right)
\\
&= 
T \left( \sum_{\alpha_\nu \in I_k^\pm} \left[ \sum_{\alpha'_{\nu'}\in I_k^\pm} \frac{1}{d_+(\alpha')d_-(\alpha')} \!\! \sum_{\alpha_\nu \overset{\beta}{\sim} \alpha'_{\nu'}}\!\! -[\beta:\alpha_\nu][\beta:\alpha'_{\nu'}] \lambda_{\alpha'_{\nu'}}\right]  {\alpha_\nu} + \cdots \Theta\right)
\\
&= 
\sum_{\alpha\in I_k} \left[  \sum_{\alpha'_{\nu'}\in I_k^\pm} \frac{1}{{d_+(\alpha')d_-(\alpha')}} \sum_{\alpha_+ \overset{\beta}{\sim} \alpha'_{\nu'}} -[\beta:\alpha_+][\beta:\alpha'_{\nu'}] \lambda_{\alpha'_{\nu'}} \right] \cdot \alpha
\\
&\quad- 
\sum_{\alpha\in I_k} \left[  \sum_{\alpha'_{\nu'}\in I_k^\pm} \frac{1}{{d_+(\alpha')d_-(\alpha')}} \sum_{\alpha_- \overset{\beta}{\sim} \alpha'_{\nu'}} -[\beta:\alpha_-][\beta:\alpha'_{\nu'}] \lambda_{\alpha'_{\nu'}} \right] \cdot \alpha.
\end{align*}
Now we use that $[\beta:\alpha] = [\beta:\alpha_+] = -[\beta:\alpha_-]$ and $[\beta:\alpha][\beta:\alpha']\neq 0$ only if either $\alpha_\pm\overset{\beta}{\sim} \alpha'_\pm$ or $\alpha_\pm \overset{\beta}{\sim} \alpha'_\mp$ together with $-[\beta:\alpha_\nu][\beta:\alpha'_{\nu'}] = \nu\nu'd(\alpha,\alpha',\beta)$ to find that 
\begin{align*}
TP\left(\sum_{s\in I_k^*} \lambda_{s}\cdot s \right)
&= \sum_{\alpha\in I_k} \left[  \sum_{\alpha'_{\nu'}\in I_k^\pm} \frac{1}{{d_+(\alpha')d_-(\alpha')}} \sum_{\alpha_+ \overset{\beta}{\sim} \alpha'_{\nu'}} d(\alpha,\alpha',\beta) \cdot \nu'\lambda_{\alpha'_{\nu'}} \right] \cdot \alpha
\\
&\quad + 
\sum_{\alpha\in I_k} \left[  \sum_{\alpha'_{\nu'}\in I_k^\pm} \frac{1}{{d_+(\alpha')d_-(\alpha')}} \sum_{\alpha_- \overset{\beta}{\sim} \alpha'_{\nu'}}  d(\alpha,\alpha',\beta) \cdot \nu'\lambda_{\alpha'_{\nu'}} \right] \cdot \alpha
\\
&=
\sum_{\alpha\in I_k} \left[  \sum_{\alpha\neq \alpha'\in I_k} \frac{1}{{d_+(\alpha')d_-(\alpha')}} \sum_{\beta\in I_{k+1}} d(\alpha,\alpha',\beta) \cdot (\lambda_{\alpha'_{+}}-\lambda_{\alpha'_{-}}) \right] \cdot \alpha
\end{align*}
showing the first equality. For the second equality, we have $TI = \Id$, hence $TPI = BTI = B$ and lastly $B^n = B^n TI = TP^n I$.
\end{proof}

\corollary{
For the operators $B$ and $P$ defined above, $n\in \N$ and all $\alpha,\alpha' \in I_k$,
$$ \langle B^n(\alpha),\alpha'\rangle = \langle P^n(\alpha_+), \alpha'_+\rangle - \langle P^n(\alpha_+), \alpha'_-\rangle. $$
}
\begin{proof}
Using $B^n = TP^n I$, we compute the coefficient $\langle B^n(\alpha),\alpha'\rangle$ of $\alpha'$ in $B(\alpha)$ by 
\begin{align*}
\langle B^n(\alpha),\alpha'\rangle 
&= \langle TP^n I (\alpha),\alpha'\rangle
= \langle TP^n (\alpha_+),\alpha'\rangle
= \langle P^n(\alpha_+) ,\alpha'_+ \rangle - \langle P^n(\alpha_+) ,\alpha'_- \rangle. \qedhere
\end{align*}
\end{proof}

In particular, we can define the following quantities generalising the idea of return probabilities.

\definition{
For the random walk described by $P$ and $\alpha\in I_k$, we define the return probabilities $p_{\alpha,+}$ and the probabilities of returning with reversed orientation $p_{\alpha,-}$ respectively by
$$ p_{\alpha,+}(n) = \langle P^n(\alpha_+), \alpha_+\rangle 
\quad \text{and} \quad 
p_{\alpha,-}(n) = \langle P^n(\alpha_+), \alpha_-\rangle. $$
For the process described by $B$ we define 
$$ p_\alpha(n) = \langle B^n(\alpha), \alpha\rangle. $$
}
Note that all three quantities are independent of the choice of preferred orientations.
Notice also that $p_\alpha(n) =  p_{\alpha,+}(n) - p_{\alpha,-}(n)$, hence, after summing over all $G$-types of $k$-cells we obtain the following expression for the von Neumann trace of $B^n$.
\corollary{
For $n\in \N$, the von Neumann trace of $B^n$ is given by 
$$ \tr_{\cal NG} (B^n) = \sum_{\alpha \in I_k(G\setminus X)} p_{\alpha}(n) = \sum_{\alpha\in I_k(G\setminus X)} p_{\alpha,+}(n) - p_{\alpha,-}(n). $$
}

Before moving on and relating this random walk to the upper $k$-Laplacian, we will introduce one extra parameter that will prove useful later on when studying the spectra of the operators.
		\subsection{Lazy Degree \texorpdfstring{$k$}{k}-Upper Random Walks} 
			Given a random walk $\frak R = (V, P)$ with state space $V$ and propagation operator $P$ and $q\in [0,1]$, the $q$-lazy version $\frak R_q$ of the random walk $\frak R$
is given by $\frak R_q = (V, q\Id + (1+q)P)$, that is the random walk on $V$ that stays put with probability $q$ and moves according to the random walk $\frak R$ with probability $1-q$. 
In particular, the moving probabilities for the $q$-lazy version $\frak R^k_q(X)$ of the previously defined random walk on $I_k^*$ are given as follows:
\begin{itemize}
\item For the absorbing state $\Theta$,
$ \qquad \P_q(\Theta\to \Theta) = 1, \qquad \P_q(\Theta\to \alpha_\nu) = 0. $
\item For $\alpha_\nu \in I_k^\pm$ (that is $\alpha\in I_k$ and $\nu\in \{\pm\}$),
\begin{align*}
&\P_q(\alpha_\nu \to \alpha_\nu) = q, \qquad
\P_q(\alpha_\nu \to -\alpha_\nu) = 0. 
\end{align*}
\item For $\alpha_\nu,\alpha'_{\nu'}\in I_k^\pm$ with $\alpha_\nu \neq \pm \alpha'_{\nu'}$,
$$\P_q(\alpha_\nu \to \alpha'_{\nu'}) = \frac{1-q}{{d_+(\alpha)d_-(\alpha)}}\sum_{\stackrel{\beta\in I_{k+1}}{\alpha_\nu \overset{\beta}{\sim} \alpha'_{\nu'}}} d(\alpha_\nu,\alpha_{\nu'},\beta) = (1-q)\P(\alpha_\nu \to \alpha'_{\nu'}). $$
\item Lastly,
$$ \P_q(\alpha_\nu \to \Theta) = 1-\sum_{\alpha'_{\nu'}\in I_k^\pm} \P_q(\alpha_\nu \to \alpha'_{\nu'}) = (1-q)\P(\alpha_\nu \to \Theta). $$
\end{itemize}

In the same spirit, we define $B_q\acts \ell^2 I_k$ by $B_q = q\Id + (1-q)B$. 
This operator is given by
\begin{align*}
B_q\left( \sum_{\alpha\in I_k} \lambda_\alpha\cdot \alpha\right) &= \sum_{\alpha\in I_k} \left[ q\lambda_\alpha + \sum_{\alpha\neq \alpha'\in I_k} \frac{1-q}{d_+(\alpha')d_-(\alpha')} \sum_{\beta\in I_{k+1}} d(\alpha,\alpha',\beta)\lambda_{\alpha'} \right] \cdot \alpha.
\end{align*}

\corollary{
These operators satisfy $B_q T = TP_q$, $B_q = TP_q I$ and $B_q^n = TP_q^nI$.
}
\begin{proof}
This follows from the previous equalities since $B_qT = qT+(1-q)BT = qT+(1-q)TP = TP_q$ and $TP_q I = qTI + (1-q)TPI = q + (1-q)B = B_q$. 
\end{proof}

As before, we consider the probabilities of returning to the same $k$-cell with the same orientation or the reversed orientation respectively.
\definition{
For $\alpha\in I_k$, we define the quantities
\begin{align*}
p_{q,\alpha,+}(n) &= \langle P_q^n(\alpha_+), \alpha_+\rangle,
\\
p_{q,\alpha,-}(n) &= \langle P_q^n(\alpha_+), \alpha_-\rangle,
\\
p_{q,\alpha}(n) &= \langle B_q^n(\alpha), \alpha\rangle  
\end{align*}
}

Again, 
$p_{q,\alpha}(n) = p_{q,\alpha,+}(n) - p_{q,\alpha,-}(n) $, hence we
can compute the von Neumann trace of $B_q^n$ using the probabilities of the random walk.
We define  
$$ p_q(n) = \tr_{\cal NG}(B_q^n) = \sum_{\alpha\in I_k(G\setminus X)} p_{q,\alpha,+}(n) - p_{q,\alpha,-}(n) $$
	\section{Relationship to Laplace Operators}
		We now compare the operator $B_q$ to the upper Laplacian $\Delta^\up_k$. 
Recall from Equation~(\ref{E_Dup_Formula}) that $\Delta^\up_k$ acts on $\ell^2 I_k$ by
\begin{align*}
\Delta^\up_k\left( \sum_{\alpha\in I_k} \lambda_\alpha\cdot \alpha\right) 
&=
\sum_{\alpha\in I_k}\left[ \sum_{\beta\in I_{k+1}} [\beta:\alpha]^2 \lambda_\alpha - \sum_{\alpha'\neq \alpha\in I_k}\sum_{\beta\in I_{k+1}} -[\beta:\alpha][\beta:\alpha'] \lambda_{\alpha'} \right] \cdot \alpha
\\
&= 
\sum_{\alpha\in I_k}\left[ d_{+,2}(\alpha) \lambda_\alpha - \sum_{\alpha'\neq \alpha\in I_k}\sum_{\beta\in I_{k+1}} d(\alpha,\alpha',\beta) \lambda_{\alpha'} \right] \cdot \alpha
\end{align*}

\theorem{
\label{Thm_Bq_Id_Delta}
Let $B_q\acts \ell^2 I_k = \ell^2C^\cell_k(X)$ be the operator $B_q = TP_q I$ defined as above. Then 
\begin{align*}
B_q \circ M_{1,q} = \Id - \Delta^\up_k \circ M_{2,q}, 
\end{align*} 
where $M_{1,q}, M_{2,q}\acts \ell^2 I_k$ are the non-negative multiplication operators given by
\begin{align*}
M_{1,q} = &\frac{d_+d_-}{qd_+d_-+(1-q)d_{+,2}},
\\
&\sum_{\alpha\in I_k} \lambda_\alpha \alpha \mapsto \sum_{\alpha\in I_k} \frac{d_+(\alpha)d_-(\alpha)}{qd_+(\alpha)d_-(\alpha) + (1-q)d_{+,2}(\alpha)}\cdot\lambda_\alpha \alpha,\\
M_{2,q} = &\frac{1-q}{qd_+d_-+(1-q)d_{+,2}},\\
&\sum_{\alpha\in I_k} \lambda_\alpha \alpha \mapsto  \sum_{\alpha\in I_k} \frac{1-q}{qd_+(\alpha)d_-(\alpha) + (1-q)d_{+,2}(\alpha)}\cdot\lambda_\alpha \alpha.
\end{align*}
}

\begin{proof}
For $\alpha, \alpha'\in I_k$ we compare the contributions $(B_q\circ M_{1,q})_{\alpha',\alpha}$ and $(\Id-\Delta^\up_k\circ M_{2,q})_{\alpha',\alpha}$ coming from the coefficient of $\alpha$ in the argument to the coefficient of $\alpha'$ in the image.\footnote{In the sense that both are operators $\Xi\acts \ell^2C_k^\cell(X)$ that can be written as $$\Xi\colon \sum_{\alpha\in I_k} \lambda_\alpha \alpha \mapsto \sum_{\alpha'\in I_k}\left[\sum_{\alpha\in I_k} \Xi_{\alpha',\alpha} \lambda_{\alpha}\right]\cdot {\alpha'}.$$}
 For $\alpha\neq \alpha'$ these contributions are given by
\begin{align*}
(B_q\circ M_{1,q})_{\alpha',\alpha} 
&= \frac{d_+(\alpha)d_-(\alpha)}{qd_+(\alpha)d_-(\alpha) + (1-q)d_{+,2}(\alpha)} \cdot (B_q)_{\alpha',\alpha} 
\\
&= \frac{d_+(\alpha)d_-(\alpha)}{qd_+(\alpha)d_-(\alpha) + (1-q)d_{+,2}(\alpha)} \cdot \frac{1-q}{d_+(\alpha)d_-(\alpha)}\sum_{\beta\in I_{k+1}} d(\alpha,\alpha',\beta)
\\
&= \frac{1-q}{qd_+(\alpha)d_-(\alpha) + (1-q)d_{+,2}(\alpha)} \sum_{\beta\in I_{k+1}} d(\alpha,\alpha',\beta)
\\
&= 0 - \frac{1-q}{qd_+(\alpha)d_-(\alpha) + (1-q)d_{+,2}(\alpha)} \cdot \left(- \sum_{\beta\in I_{k+1}} d(\alpha,\alpha',\beta) \right)
\\
&= (\Id-\Delta^\up_k\circ M_{2,q})_{\alpha',\alpha}
\end{align*}
and for $\alpha=\alpha'$ by
\begin{align*}
(B_q\circ M_{1,q})_{\alpha,\alpha} 
&= 
 q\cdot \frac{d_+(\alpha)d_-(\alpha)}{qd_+(\alpha)d_-(\alpha) + (1-q)d_{+,2}(\alpha)} 
 \\
 &=
 1 - \frac{(1-q)d_{+,2}(\alpha)}{qd_+(\alpha)d_-(\alpha) + (1-q)d_{+,2}(\alpha)}
 = (\Id - \Delta^\up_k\circ M_{2,q})_{\alpha,\alpha}. 
\end{align*}
Since all these coefficients agree, the claim follows.
\end{proof}

\remark{
\label{R_simplicial_upper_k_regular}
The construction here generalises the one given by O.~Parzanchevski and R.~Rosenthal on simplicial complexes in~\cite{PaRo13} and the previous theorem generalises Proposition~2.8~(1) of their paper. Considering the random walk of O.~Parzanchevski and R.~Rosenthal in degree $k=(d-1)$, the incidence numbers of a simplicial complex (viewed as a CW complex) are given as $[\beta:\alpha]\in\{0,\pm 1\}$, where $\pm 1$ occurs if the $(d-1)$-simplex $\alpha$ is in the boundary of the $d$-simplex $\beta$, with sign depending on orientations. Therefore, 
$$ d_+(\alpha) = d_{+,2}(\alpha) = \deg(\alpha), \quad d_-(\alpha) = d, $$
where $\deg(\alpha)$ denoted the number of $d$-simplices $\beta\in I_d$ containing $\alpha$. 
Hence, 
$$ \frac{d}{q(d-1)+1}B_q = \Id - \frac{1-q}{q(d-1)+1} \cdot \frac{\Delta^\up_{d-1}}{\deg(\alpha)}, $$
where $\Delta^\up_w = \frac{\Delta^\up_{d-1}}{\deg(\alpha)}$ is the weighted upper Laplacian used by O.~Parzanchevski and R.~Rosenthal, defined by using a weighted scalar product on $\ell^2 I_k$.
Note that in this case, the diagonal operator $M_{1,q}$ is given by multiplication by a constant (depending on $q$ and $d$ but not on $\alpha\in I_{d-1}$). 
}
	\section{Application to \texorpdfstring{$L^2$}{L2}-Invariants}
		As a short reminder on $L^2$-invariants, we briefly recall some important definitions. 
For a proper introduction we refer to \cite{KamBook} and \cite{LckL2}.

Let $M$ be a manifold with cocompact free and proper action by a (discrete) group $G$. 
We denote by $\cal NG$ the von Neumann algebra of $G$. 
On $\ell^2 G$, we define the von Neumann trace by $\tr_{\cal NG}(f) = \langle f(e), e\rangle_{\ell^2 G}$.
This notion can be extended to Hilbert $\cal NG$-modules, that is, Hilbert spaces $V$ with isometric left $G$-action that have an isometric $G$-embedding $V\hookrightarrow H\otimes \ell^2 G$ for some Hilbert space $H$. 

\definition{
Given two Hilbert $\cal NG$-modules $U, V$ and a closed densely defined
$G$-equivariant operator $f\colon U\to V$, we define the spectral density function $F(f)$ of $f$ using the family of $G$-equivariant spectral projectors of the self-adjoint operator $f^*f$ $\left\{E^{f^*f}_\lambda\right\}_{\lambda\geq 0}$ by
$$ F(f)\colon \R_{\geq 0} \to [0,\infty], \qquad F(f)(\lambda) = \tr_{\cal NG}\left( E^{f^*f}_{\lambda^2} \right).$$
We call $f$ Fredholm if there is $\lambda>0$ with $F(f)(\lambda) <\infty$. 
}

This function measures the content of the spectrum of the operator $f$ close to zero in terms of the parameter $\lambda$. 
Notice that $F$ is monotonously increasing and right-continuous.
The value at zero, $F(f)(0)$, measures the size of the kernel of $f$ and is called the $L^2$-Betti number. 
Additionally, we define the Novikov-Shubin invariants to measure the spectral content close to zero.

\definition{
	Let $F$ be monotonously increasing and right-continuous. 
	We define the $L^2$-Betti number $b^{(2)}(F) = F(0)$ and the Novikov-Shubin invariant 
	$$ \alpha(F) = \liminf_{\lambda\searrow 0} \frac{\log(F(\lambda) - F(0))}{\log(\lambda)}. $$
	If $F(f)$ is the spectral density function of $f$, we also write $b^{(2)}(f)=b^{(2)}(F(f))$ and $\alpha(f) = \alpha(F(f))$.
}

Here, $\alpha(F)$ measures the asymptotic behaviour of $F$ close to zero. 
Indeed, if asymptotically $F(\lambda) - F(0) \sim \lambda^{\alpha}$ as $\lambda\searrow 0$, then $\alpha=\alpha(F)$.
In particular, we only care for the asymptotic behaviour of $F$ close to zero. 
This motivates the following equivalence relation defined on such functions. 

\definition{
Let $F, F'\colon \R_{\geq 0}\to [0,\infty]$ be two monotonously increasing right-continuous functions. 
We call $F$ and $F'$ dilatationally equivalent, and write $F\sim F'$, if there are constants $C > 0$ and $\lambda_0 > 0$ such
that
$$ F'(C^{-1}\lambda) \leq F(\lambda) \leq F'(C\lambda)  \qquad \text{ for all } \lambda\in [0,\lambda_0].$$
}

Novikov-Shubin invariants are invariant under dilatational equivalence. 
Given a free $G$-CW complex of finite type, we use the differentials of the cellular $L^2$-chain complex, $d_k\colon \ell^2C^{\cell}_k(X;\cal NG)\to \ell^2 C^{\cell}_{k-1}(X;\cal NG)$, and the cellular Laplace operators, $\Delta_k = d_kd_k^* + d_k^*d_k \acts \ell^2C^{\cell}_k(X;\cal NG)$, to associate $L^2$-Betti numbers and Novikov-Shubin invariants to $X$. 

\definition{
Let $X$ be a free $G$-CW complex of finite type. We define the $L^2$-Betti numbers of $X$ by
$ b^{(2)}_k(X) =  b^{(2)}(\Delta_k) $
and the Novikov-Shubin invariants by 
$ \alpha_k(X) = \alpha(d_k). $
}
Here, we use $d_k$ instead of $\Delta_k$ in the definition of the Novikov-Shubin invariants, as they provide a more detailed picture. 
Indeed, $\alpha(\Delta_k) = \half \cdot \min\{\alpha_k(X), \alpha_{k+1}(X)\}$. 

		\subsection{Computing Novikov-Shubin Invariants}
			We now study the connection between the random walk $\frak R^k_q=\frak R^k_q(X)$ and the Novikov-Shubin invariant $\alpha_k(X)$ for a free $G$-CW complex $X$ of finite type. 
In degree zero it is reasonable to consider only connected spaces (since the $\ell^2$-spaces and the Laplace operator split as a direct sum with one summand for each connected component).
By the same reasoning, we can assume without loss of generality the following analogue in degree $k$.
\definition{
Let $X$ be a CW complex.
We call $X$ upper $k$-connected if $|I_k|\geq 2$ and for all $\alpha,\alpha'\in I_k$ there are 
$$\alpha=\alpha_0,\alpha_1,\dots, \alpha_{n-1},\alpha_n=\alpha'\in I_k \quad\text{ and }\quad \beta_1,\dots, \beta_{n}\in I_{k+1}$$
 such that  $[\beta_i,\alpha_{i-1}]\neq 0$ and $[\beta_i:\alpha_i]\neq 0$ for all $1\leq i\leq n$, that is $\beta_i$ is attached non-trivially to $\alpha_{i-1}$ and $\alpha_i$.
} 
This condition implies for $\frak R^k_q$ that a random walker can move from any $k$-cell $\alpha_\pm\in I^\pm_k$ to any other (unoriented) $k$-cell $\alpha'$  (that is, to one of the oriented $k$-cells $\alpha'_+$ or $\alpha'_-$).
Furthermore, we get bounds on the quantities from Definition~\ref{Def_RW_ds}.
\lemma{
Let $X$ be an upper $k$-connected free $G$-CW complex of finite type. 
Then there exists $D\geq 1$ such that 
$$D\geq d_{+,2}, d_+, d_-\geq 1.$$
In particular, if $q\in [0,1)$ then the operators $M_{1,q}$ and $M_{2,q}$ are positive multiplication operators bounded from below by 
$$M_{1,q} \geq D^{-2}>0 \quad \text{and} \quad M_{2,q} = (1-q)D^{-2}>0.$$
}

\begin{proof}
Let $\alpha\in I_k$ be arbitrary and let $\alpha\neq \alpha'\in I_k$ be any other $k$-cell.
Since $X$ is upper $k$-connected, by definition there is a sequence of $k$-cells $\alpha_i\in I_k$ and $(k+1)$-cells $\beta_i\in I_{k+1}$ connecting $\alpha$ to $\alpha'$.
In particular, there exists a $(k+1)$-cell $\beta_1\in I_{k+1}$ such that $[\beta_1:\alpha]\neq 0$ and a $k$-cell $\alpha\neq \alpha_1\in I_k$ such that $[\beta_1:\alpha_1]\neq 0$. 
Hence
\begin{align*}
d_{+,2}(\alpha) &\geq [\beta_1:\alpha]^2 \geq 1, \\
d_+(\alpha) &\geq |[\beta_1:\alpha]| \geq 1, \\
d_-(\alpha) &\geq d_-(\beta_1; \alpha) \geq |[\beta_1:\alpha_1]|\geq 1. 
\end{align*}
Since $X$ is of finite type and these quantities depend only on the $G$-type of $\alpha$, there exists 
$$ D = \sup_{\alpha\in I_k}  \{d_{+,2}(\alpha), d_{+}(\alpha), d_{-}(\alpha) \} = \max_{\alpha\in I_k(G\setminus X)} \{d_{+,2}(\alpha), d_{+}(\alpha), d_{-}(\alpha) \} \geq 1. $$
It follows, therefore, that
\begin{align*}
M_{1,q} &= \frac{d_+ d_-}{qd_+ d_- + (1-q) d_{+,2}} \geq \frac{1}{qD^2 + (1-q)D} \geq \frac{1}{D^2}>0, \\
M_{2,q} &= \frac{1-q}{qd_+ d_- + (1-q) d_{+,2}} \geq \frac{1-q}{qD^2 + (1-q)D} \geq \frac{1-q}{D^2}>0,
\end{align*}
and the claim follows.
\end{proof}

Generalising the notion of regular graphs, we introduce the following notion of upper $k$-regular free $G$-CW complexes.
\definition{
Let $X$ be a free $G$-CW complex of finite type.
We call $X$ upper $k$-regular if $X$ is upper $k$-connected and $d_+d_- = d_+(\alpha)d_-(\alpha)$ and $d_{+,2} = d_{+,2}(\alpha)$ are independent of the cell $\alpha\in I_k$.
}

In this case, also the multiplication operators $M_{1,q}$ and $M_{2,q}$ are just multiplication with a constant. Hence, the formula connecting $B_q$ and $\Delta^\up_k$ simplifies further.
\example{
This happens, for example, in the case where $X$ is a regular graph for $k=0$ or in the case studied in \cite{BCDS23} for $k=\dim(X)-1$. 
While we can compute the Novikov-Shubin invariants in the latter case by Poincaré duality, we can use their values to obtain informations about the random walk considered.
}

\corollary{
Let $X$ be an upper $k$-regular $G$-CW complex of finite type and $q\in [0,1]$. 
Then
$$ C_{1,q} B_q = \Id - C_{2,q} \Delta^\up_k, $$
where $C_{1,q}$ and $C_{2,q}$ are positive constants given as
$$C_{1,q} = \frac{d_+d_-}{qd_+d_- + (1-q)d_{+,2}} > 0 \qquad \text{and} \qquad C_{2,q} = \frac{1-q}{qd_+d_- + (1-q)d_{+,2}} > 0.$$
}

\definition{
Let $X$ be an upper $k$-regular free $G$-CW complex of finite type. We define 
\begin{align*}
\widetilde{B_q} &= C_{1,q} B_q \quad \text{ and } \quad
\widetilde{\Delta^\up_{q,k}} =  C_{2,q} \Delta^\up_k
\end{align*}
so that we have the equality 
$\widetilde{B_q} = \Id-\widetilde{\Delta^\up_{q,k}}.$
}

We now can derive bounds on the spectrum $\sigma\left(\widetilde{\Delta^\up_{q,k}}\right)$ of the operator $\widetilde{\Delta^\up_{q,k}}$. 
It is well-known that $\Dup_k$ is bounded, 

\lemma{\label{L:DeltaBounded}
Let $X$ be a free $G$-CW complex of finite type, then $\Dup_k\acts \ell^2C^\cell_k(X)$ is bounded and $\sigma(\Dup_k) \subset [0, S_{k}]$ for some $S_k<\infty$. 
}

While the precise bounds won't matter for this paper, we can find explicit bounds by straight-forward computations, 
for example, 
$$ S_{k} = \max_{\alpha\in I_k(G\backslash X)} \left\{  \sum_{\beta\in I_{k+1}} \sum_{\alpha'\in I_k}   |d(\alpha,\alpha',\beta)| \right\} < \infty.  $$  
Using this, we can prove the following lemma.

\lemma{
\label{L_DeltaTildeSpectrum01}
Let $X$ be an upper $k$-regular free $G$-CW complex of finite type. Then there exists $q_0\in (0,1)$ such that for all $q_0\leq q\leq 1$ the spectrum of $\widetilde{\Delta^\up_{q,k}}$ is contained in the unit interval,
$\sigma(\widetilde{\Delta^\up_{q,k}}) \subset [0,1]. $
}
\begin{proof}
By Lemma~\ref{L:DeltaBounded}, $\sigma(\Delta^\up_k) \subset [0,S]$ for some $S>0$, hence $\sigma(\widetilde{\Delta^\up_{q,k}}) \subset [0, C_{2,q}S]$. Note that $d_+d_-\ge 1$ and $d_{+,2}\geq 1$ so $C_{2,q} = \frac{1-q}{qd_+d_- + (1-q)d_{+,2}}$ is continuous in $q\in (0,1)$ and converges to $0$ as $q\nearrow 1$. 
In particular, there is $q_0\in (0,1)$ such that $0< C_{2,q} \leq S^{-1}$ for all $q_0\leq q\leq 1$.
\end{proof}

\corollary{
Let $X$ be an upper $k$-regular free $G$-CW complex of finite type and $q\in [q_0, 1)$.
Let $\widetilde d_{k+1} = \sqrt{C_{2,q}} d_{k+1}$.
Then $\widetilde d_{k+1}^* = \sqrt{C_{2,q}} d_{k+1}^*$ and 
$$ \widetilde{\Delta^\up_{q,k}} = \widetilde d_{k+1}\widetilde d_{k+1}^* $$
is a self-adjoint positive operator with $\sigma(\widetilde{\Delta^\up_{q,k}}) \subset [0,1]$.
}

\remark{
Recall that for $q\in [q_0,1)$, since $d$ and $\widetilde{d}$ differ only by a constant factor $\sqrt{C_{2,q}}$, their spectral density functions are dilatationally equivalent and hence their Novikov-Shubin invariants agree, that is,
$\alpha_k(X) = \alpha(d_{k+1}) = \alpha(\widetilde d_{k+1}) = \alpha({\widetilde d_{k+1}}^*)$.
}

\lemma{
\label{L_trBq_F}
Let $\chi_I$ denote the indicator function of the interval $I$, then
$$ \tr_{\cal NG} (\chi_{[1-\lambda,1]}(\widetilde{B_q})) = F(\widetilde d_{k+1}^*)(\sqrt \lambda).$$
}
\begin{proof}
Recall that $\widetilde{B_q} = \Id-\widetilde d_{k+1} \widetilde{d}_{k+1}^*$, hence
\begin{align*}
\tr_{\cal NG} \left(\chi_{[1-\lambda,1]}\left(\widetilde{B_q}\right)\right) 
&= \tr_{\cal NG}\left(\chi_{[0,\lambda]}\left(\widetilde d_{k+1} \widetilde{d}_{k+1}^*\right)\right) 
\\
&= \tr_{\cal NG}\left( E_\lambda^{\widetilde d_{k+1} \widetilde{d}_{k+1}^*}\right) 
= F\left(\widetilde{d}_{k+1}^*\,\right)\left(\sqrt{\lambda}\right). \qedhere
\end{align*}
\end{proof}

We can now proceed in the same way as in degree zero, compare Lück's book~\cite[§2.1.4]{LckL2}.

\theorem{
\label{Thm_RW_ak}
Let $X$ be an upper $k$-regular free $G$-CW complex of finite type and $q\in [q_0, 1)$, with $q_0$ given by Lemma~\ref{L_DeltaTildeSpectrum01}.
 Then $\alpha_k(X) = 2a$ if and only if there is a constant $C>0$ such that for all $n\in \N$,
$$  C_{1,q}^{-n} \left( b^{(2)}(d_{k+1}^*) + C^{-1}n^{-a}\right) \quad \leq \quad p_q(n) \quad\leq\quad  C_{1,q}^{-n}\left(b^{(2)}(d_{k+1}^*) + Cn^{-a}\right). $$
}
\begin{proof}
Since by Lemma~\ref{L_DeltaTildeSpectrum01}, $\sigma(\widetilde{\Delta^\up_{q,k}}) \subset [0,1]$ 
and by construction $ \widetilde{\Delta^\up_{q,k}} = \Id - \widetilde{B_q}$, it follows that also
$\sigma(\widetilde{B_q})\subset [0,1]$.
Therefore,
$$ (1-\lambda)^n\chi_{[1-\lambda,1]}(\widetilde{B_q}) 
\quad \leq \quad 
\widetilde{B_q}^n 
\quad \leq \quad 
(1-\lambda)^n\chi_{[0,1-\lambda]}(\widetilde{B_q}) + \chi_{[1-\lambda,1]}(\widetilde{B_q}). $$
Taking traces using Lemma~\ref{L_trBq_F} and denoting $\widetilde{p}_q(n) = \tr_{\cal NG}(\widetilde{B_q}^n)$ yields
$$
(1-\lambda)^n F(\widetilde d_{k+1}^*)(\sqrt \lambda) 
\quad \leq \quad 
\widetilde{p}_q(n)
\quad \leq \quad 
(1-\lambda)^n + F(\widetilde d_{k+1}^*)(\sqrt \lambda).
$$ 
By rearranging these terms and taking logarithms, we obtain the inequalities
\begin{align}
\label{Eq_LckRWEqs_1} 
\frac{\log\left( F(\widetilde d_{k+1}^*)(\sqrt \lambda)- b^{(2)}(d_{k+1}^*)\right) }{\log \lambda} 
&\leq \frac{\log(\widetilde{p}_q(n)- (1-\lambda)^n b^{(2)}(d_{k+1}^*))}{\log \lambda} - n\cdot  \frac{\log(1-\lambda)}{\log \lambda}, \\ 
\label{Eq_LckRWEqs_2}
\frac{\log\left( F(\widetilde d_{k+1}^*)(\sqrt \lambda)- b^{(2)}(d_{k+1}^*)\right)}{\log \lambda}  
&\geq \frac{\log(\widetilde{p}_q(n)- b^{(2)}(d_{k+1}^*) - (1-\lambda)^n)}{\log \lambda}. 
\end{align}
Using $b^{(2)}(d_{k+1}^*) = b^{(2)}(\widetilde d_{k+1}^*)$ and taking the limit inferior for $\lambda\searrow 0$ on the left-hand-sides yields
\begin{align*}
\liminf_{\lambda\searrow 0} \frac{\log\left( F(\widetilde d)(\sqrt \lambda)- b^{(2)}(d_{k+1}^*)\right)}{\log \lambda} 
&= \frac{\alpha(\widetilde d_{k+1}^*)}{2} =
\frac{\alpha(\widetilde d_{k+1})}{2} =
 \frac{\alpha_k(X)}{2}.
\end{align*} 
After substituting $p(n) = \widetilde{p}_q(n)- b^{(2)}(d_{k+1}^*)$, the term on the right-hand-side of Equation~(\ref{Eq_LckRWEqs_2}) agrees with the term in~\cite[Thm. 2.48]{LckL2}, so by the same argument 
$$ \alpha_k(X) \leq 2a \qquad \text{if } \widetilde{p}_q(n)\geq b^{(2)} (d_{k+1}^*) + Dn^{-a} \text{ for } n\geq 1,$$
for some constant $D>0$.

For the right-hand-side of Equation~(\ref{Eq_LckRWEqs_1}), let $\varepsilon>0$ be arbitrarily small and $n=n(\lambda)$ the largest integer such that $n\leq \lambda^{-\varepsilon}$, that is $n=\lfloor \lambda^{-\varepsilon}\rfloor$. 
If $\widetilde{p}(n)\geq Cn^{-a} + b^{(2)}(d)$ for some constant $C>0$ and $n\geq 1$, we obtain
\begin{align*}
&\frac{\log(\widetilde{p}_q(n)- (1-\lambda)^n b^{(2)}(d_{k+1}^*))}{\log \lambda} - n\cdot  \frac{\log(1-\lambda)}{\log \lambda} \\
&\qquad\geq 
\frac{\log(Cn^{-a} + \left[1- (1-\lambda)^n\right] b^{(2)}(d_{k+1}^*))}{\log \lambda} - \frac{\log(1-\lambda)}{\lambda^\varepsilon \log \lambda}
\\
&\qquad\geq \frac{\log(Cn^{-a})}{\log \lambda} - \frac{\log(1-\lambda)}{\lambda^\varepsilon \log \lambda},
\end{align*}
where we use $\left[1- (1-\lambda)^n\right] b^{(2)}(d_{k+1}^*) \geq 0$ (indeed, even $\left[1- (1-\lambda)^n\right] b^{(2)}(d_{k+1}^*) \xrightarrow{\lambda\searrow 0} 1-e^{-\varepsilon} $). 
From here, we proceed precisely as in~\cite[Thm.~2.48]{LckL2} and find 
$$\alpha_k(X) \geq 2a \qquad \text{if } \widetilde{p}_q(n)\leq b^{2}(d_{k+1}^*) + Cn^{-a}\text{ for } n\geq 1, $$
concluding the proof of the theorem.
\end{proof}

\remark{
This generalises the theorem in degree zero, since in degree zero we have $d_+=d_{+,2}=|S|$, where $|S|$ is the size of a finite generating set of $G$ chosen in the construction of $\Cayley(G)$, and $d_-=1$. Thus $C_{1,q} = 1$ and the exponential decay factor $C_{1,q}^{-n} = 1$ disappears.  
}

\example{
Let $k\geq 2$ and let $G$ be a finitely generated group with Cayley graph $\Cayley(G)$. 
Construct a $G$-CW complex $X$ in the following way.
\begin{itemize}
\item Start with $X^{(1)} = \Cayley(G)$.\footnote{If $k\geq 3$ and $G$ is finitely presented, we can further glue in $2$-cells according to the relations in $G$, so that $X^{(2)}$ is the Cayley complex of $G$. In that case the constructed CW complex $X$ satisfies $\pi_1(X)=G$, see for example A. Hatcher's book~\cite[p. 77]{Hatcher}.}
\item For every $g\in G$ glue a $k$-cell $\alpha_g$ to $X^{(1)}$ by collapsing the boundary of $\alpha_g$ to the vertex $v_g$ corresponding to $g\in G$ in $\Cayley(G)$. This defines $X^{(k)}$.
\item For every edge $(g,gs)$ in the Cayley graph, glue one $(k+1)$-cell $\beta_{g,gs}$ to $X^{(k)}$ by sending the boundary of $\beta_{g,gs}$ to $\alpha_g\cup (g,gs) \cup \alpha_{gs}$ such that $[\beta_{g,gs}:\alpha_g] = -[\beta_{g,gs}:\alpha_{gs}] \in \{\pm 1\}$.
This defines $X^{(k+1)} = X$.
\end{itemize}
On $X$, the degree $k$-upper random walk $\frak R^k$ agrees with the random walk $\frak R$ on $\Cayley(G)$ when identifying the state corresponding to $\alpha_g = (\alpha_g)_+$ in $\frak R^k$ with the state corresponding to $g$ in $\frak R$. In particular, for $\frak R^k$ we have $p_-(n) \equiv 0$ so that $p(n) = p_+(n)$ is the usual return probability. 
Further, the values $d_+=d_{2,+}=|S|$ and $d_-=1$ agree with the values on $\Cayley(G)$, so that $C_{1,q}=1$. 
Therefore, the previous theorem tells us that for $X$ we obtain
$\alpha_k(X) = \alpha_0(X)$. 

Indeed, we can also see this in a different way because $d_{k+1}d_{k+1}^*$ and $d_1d_1^*$ are unitarely equivalent by identifying $\alpha_g$ with $g$ and $\beta_{g,gs}$ with $(g,gs)$.
}
		\subsection{Example: Degree 1-Upper Random Walk on \texorpdfstring{$\R^2$}{R2}}
			\begin{minipage}{0.65\textwidth}
Consider $\R^2$ as a $\Z^2$-CW complex of finite type as shown on the right, with arrows indicating the chosen preferred orientation.
For notation's sake, we will write $\Z^2$ as a multiplicative group with unit element $1\in \Z^2$. 
Let $x$ and $y$ be two generators of $\Z^2=\langle x,y \:|\: [x,y]=1\rangle$ and the $\Z^2$-action on this CW complex be generated by $x$ shifting to the right by one and $y$ shifting up by one. 
The red cells indicate $\Z^2$-bases. 
We will denote the $0$-basis $\cal B_0 = \{\gamma_\bullet\}$, the $1$-basis $\cal B_1 = \{\alpha_\uparrow, \alpha_\rightarrow\}$ and the $2$-basis $\cal B_2 = \{\beta_{\circlearrowright}\}$ in the way suggested by the indices. 
Given a cell $c$ and $g=x^ay^b\in \Z^2$, we denote by $gc$ the cell obtained by translating $c$ by $g$, that is $a$ units to the right and $b$ units up. \newline
\vspace*{-.75\baselineskip}
\end{minipage} \hfill
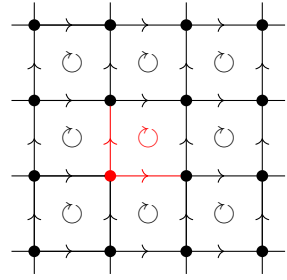
\begin{minipage}{0.35\textwidth}
\centering
\begin{tikzpicture}
\foreach \x in {-1,0,1,2}{
	\foreach \y in {-1,0,1,2}{
		\def\clr{black}
		\ifnum \x=0 \ifnum \y=0 \def\clr{red} 
		\else \fi \else \fi
		\node[fill=\clr, circle, minimum size=0.4em, inner sep=0pt, outer sep=0pt] at (\x, \y){};
		\ifnum \x=2 {} \else \ifnum \y=2 {} \else
			\node[\clr] at (\x+0.5, \y+0.5) {$\circlearrowright$};
		\fi \fi
	}
}
\foreach \x in {-1,0,1,2}{
	\foreach \y in {-1,0,1,2}{
		\def\clr{black}
		\ifnum \x=0 \ifnum \y=0 \def\clr{red} 
		\else \fi \else \fi
		\ifnum \x=2 
			\draw (-1.3, \y) -- (-0.07, \y);
			\draw (\x, \y) -- (\x+0.3, \y);		
		\else
			\draw[\clr, decoration={markings, mark=at position 0.5 with {\arrow{>}}},
        postaction={decorate}] (\x+0.07,\y) -- (\x+0.93,\y);
		\fi 
		\ifnum \y=2 
			\draw (\x,-1.3) -- (\x, -0.07);
			\draw (\x, \y+0.07) -- (\x, \y+0.3);		
		\else 
			\draw[\clr, decoration={markings, mark=at position 0.5 with {\arrow{>}}},
        postaction={decorate}] (\x, \y+0.07) -- (\x, \y+0.93);
		\fi
	}
}
\end{tikzpicture}
\end{minipage} 

The incidence numbers between a $2$-cell $\beta$ and a $1$-cell $\alpha$ are given by $[\beta:\alpha] = 0$ if $\beta$ and $\alpha$ do not touch and $[\beta:\alpha] = \pm 1$ if the cells touch; with sign $+1$ if the orientation $\beta$ induces on $\alpha$ agrees with the orientation on $\alpha$ and $-1$ otherwise.
This is an upper $2$-regular CW complex with
\begin{align*}
d_{+} = 2, \ d_{+,2} = 2, \ d_- = 3,  \quad C_{1,q} = \frac{3}{2q+1}, \quad  C_{2,q} = \frac{1}{2}\frac{1-q}{2q+1}, \quad C_{1,q}^{-1}C_{2,q} = \frac{1-q}{6}.
\end{align*}

The upper Laplacian $\Delta = \Delta^\up_1\acts \ell^2 \left((\R^2)^{(1)}\right)$ in degree one can be written, with respect to the basis $\cal B_1$, as the $\C[\Z^2]$-valued matrix
$$ \Delta = 2- \begin{pmatrix}
x+x^{-1} & 1-x-y^{-1}+xy^{-1} \\
1-x^{-1}-y+x^{-1}y & y+y^{-1}
\end{pmatrix}.
$$

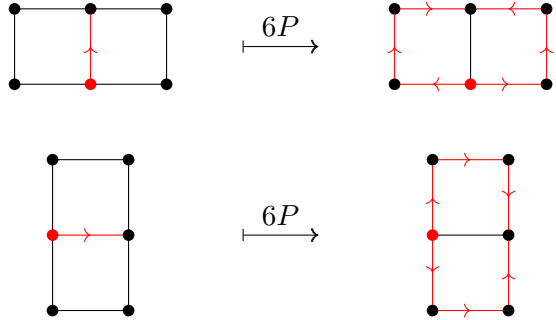
\begin{figure}[h!]
\centering
\begin{tikzpicture}
\foreach \x in {-1,1}{
	\foreach \y in {0,1}{
		\node[fill=black, circle, minimum size=0.4em, inner sep=0pt, outer sep=0pt] at (\x,\y){};
	}
}
\node[fill=black, circle, minimum size=0.4em, inner sep=0pt, outer sep=0pt] at (0,1){};
\draw (-1,0) -- (1,0) -- (1,1) -- (-1,1) -- (-1,0);
\node[fill=red, circle, minimum size=0.4em, inner sep=0pt, outer sep=0pt] at (0,0){};
\draw[red, decoration={markings, mark=at position 0.5 with {\arrow{>}}}, postaction={decorate}] (0,0.07) -- (0,0.93);
\draw[|->, decoration={markings,mark=at position 1 with {\arrow[scale=1]{>}}},
    postaction={decorate}, 
    shorten >=0.4pt] (2,0.5)--node [midway, above]{$6P$} (3,0.5);
\begin{scope}[shift={(5,0)}]
\foreach \x in {-1,1}{
	\foreach \y in {0,1}{
		\node[fill=black, circle, minimum size=0.4em, inner sep=0pt, outer sep=0pt] at (\x,\y){};
	}
}
\node[fill=black, circle, minimum size=0.4em, inner sep=0pt, outer sep=0pt] at (0,1){};
\draw (0,0)--(0,1);
\node[fill=red, circle, minimum size=0.4em, inner sep=0pt, outer sep=0pt] at (0,0){};
\draw[red, decoration={markings, mark=at position 0.5 with {\arrow{>}}}, postaction={decorate}] (-1,0.07) -- (-1,0.93);
\draw[red, decoration={markings, mark=at position 0.5 with {\arrow{>}}}, postaction={decorate}] (1,0.07) -- (1,0.93);
\draw[red, decoration={markings, mark=at position 0.5 with {\arrow{>}}}, postaction={decorate}] (-0.07,0) -- (-0.93,0);
\draw[red, decoration={markings, mark=at position 0.5 with {\arrow{>}}}, postaction={decorate}] (0.07,0) -- (0.93,0);
\draw[red, decoration={markings, mark=at position 0.5 with {\arrow{>}}}, postaction={decorate}] (0.93,1) -- (0.07,1);
\draw[red, decoration={markings, mark=at position 0.5 with {\arrow{>}}}, postaction={decorate}] (-0.93,1) -- (-0.07,1);
\end{scope}	
\begin{scope}[shift={(-0.5,-2)}]
\foreach \x in {0,1}{
	\foreach \y in {-1,1}{
		\node[fill=black, circle, minimum size=0.4em, inner sep=0pt, outer sep=0pt] at (\x,\y){};
	}
}
\node[fill=black, circle, minimum size=0.4em, inner sep=0pt, outer sep=0pt] at (1,0){};
\draw (0,-1) -- (0,1) -- (1,1) -- (1,-1) -- (0,-1);
\node[fill=red, circle, minimum size=0.4em, inner sep=0pt, outer sep=0pt] at (0,0){};
\draw[red, decoration={markings, mark=at position 0.5 with {\arrow{>}}}, postaction={decorate}] (0.07,0) -- (0.93,0);
\draw[|->, decoration={markings,mark=at position 1 with {\arrow[scale=1]{>}}},
    postaction={decorate}, 
    shorten >=0.4pt] (2.5,0)--node [midway, above]{$6P$} (3.5,0);
\begin{scope}[shift={(5,0)}]
\foreach \x in {0,1}{
	\foreach \y in {-1,1}{
		\node[fill=black, circle, minimum size=0.4em, inner sep=0pt, outer sep=0pt] at (\x,\y){};
	}
}
\node[fill=black, circle, minimum size=0.4em, inner sep=0pt, outer sep=0pt] at (1,0){};
\draw (0,0)--(1,0);
\node[fill=red, circle, minimum size=0.4em, inner sep=0pt, outer sep=0pt] at (0,0){};
\draw[red, decoration={markings, mark=at position 0.5 with {\arrow{>}}}, postaction={decorate}] (0.07,-1) -- (0.93,-1);
\draw[red, decoration={markings, mark=at position 0.5 with {\arrow{>}}}, postaction={decorate}] (0.07,1) -- (0.93,1);
\draw[red, decoration={markings, mark=at position 0.5 with {\arrow{>}}}, postaction={decorate}] (0,-0.07) -- (0,-0.93);
\draw[red, decoration={markings, mark=at position 0.5 with {\arrow{>}}}, postaction={decorate}] (0,0.07) -- (0,0.93);
\draw[red, decoration={markings, mark=at position 0.5 with {\arrow{>}}}, postaction={decorate}] (1,0.93) -- (1,0.07);
\draw[red, decoration={markings, mark=at position 0.5 with {\arrow{>}}}, postaction={decorate}] (1,-0.93) -- (1,-0.07);
\end{scope}	
\end{scope}	
\end{tikzpicture}
\caption{Visual representation of propagation operator $P$}
\label{F:R2Z2CW_Moves}
\end{figure}

For the non-lazy random walk on $1$-cells, on $\cal B_1$ the propagation operator is given as described in Figure~\ref{F:R2Z2CW_Moves}. 

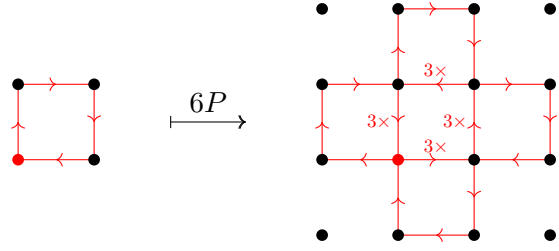
\begin{figure}[h!]
\centering
\begin{tikzpicture}
\draw[red, decoration={markings, mark=at position 0.5 with {\arrow{>}}}, postaction={decorate}] (0,0)--(0,1);
\draw[red, decoration={markings, mark=at position 0.5 with {\arrow{>}}}, postaction={decorate}] (0,1)--(1,1);
\draw[red, decoration={markings, mark=at position 0.5 with {\arrow{>}}}, postaction={decorate}] (1,1)--(1,0);
\draw[red, decoration={markings, mark=at position 0.5 with {\arrow{>}}}, postaction={decorate}] (1,0)--(0,0);
\foreach \x in {0,1}{
	\foreach \y in {0,1}{
		\def\clr{black}
		\ifnum \x=0 \ifnum \y=0 \def\clr{red} 
		\else \fi \else \fi
		\node[fill=\clr, circle, minimum size=0.4em, inner sep=0pt, outer sep=0pt] at (\x,\y){};
	}
}
\draw[|->, decoration={markings,mark=at position 1 with {\arrow[scale=1]{>}}},
    postaction={decorate}, 
    shorten >=0.4pt] (2,0.5)--node [midway, above]{$6P$} (3,0.5);
\begin{scope}[shift={(5,0)}]
\draw[red, decoration={markings, mark=at position 0.5 with {\arrow{>}}}, postaction={decorate}] (0,0) -- (-1,0);
\draw[red, decoration={markings, mark=at position 0.5 with {\arrow{>}}}, postaction={decorate}] (-1,0) -- (-1,1);
\draw[red, decoration={markings, mark=at position 0.5 with {\arrow{>}}}, postaction={decorate}] (-1,1) -- (0,1);
\draw[red, decoration={markings, mark=at position 0.5 with {\arrow{>}}}, postaction={decorate}] (0,1) -- (0,2);
\draw[red, decoration={markings, mark=at position 0.5 with {\arrow{>}}}, postaction={decorate}] (0,2) -- (1,2);
\draw[red, decoration={markings, mark=at position 0.5 with {\arrow{>}}}, postaction={decorate}] (1,2) -- (1,1);
\draw[red, decoration={markings, mark=at position 0.5 with {\arrow{>}}}, postaction={decorate}] (1,1) -- (2,1);
\draw[red, decoration={markings, mark=at position 0.5 with {\arrow{>}}}, postaction={decorate}] (2,1) -- (2,0);
\draw[red, decoration={markings, mark=at position 0.5 with {\arrow{>}}}, postaction={decorate}] (2,0) -- (1,0);
\draw[red, decoration={markings, mark=at position 0.5 with {\arrow{>}}}, postaction={decorate}] (1,0) -- (1,-1);
\draw[red, decoration={markings, mark=at position 0.5 with {\arrow{>}}}, postaction={decorate}] (1,-1) -- (0,-1);
\draw[red, decoration={markings, mark=at position 0.5 with {\arrow{>}}}, postaction={decorate}] (0,-1) -- (0,0);
\draw[red, decoration={markings, mark=at position 0.5 with {\arrow{>}}}, postaction={decorate}] (0,1) --node[midway, left, red, scale=0.65]{$3\times$} (0,0);
\draw[red, decoration={markings, mark=at position 0.5 with {\arrow{>}}}, postaction={decorate}] (1,1) --node[midway, above, red, scale=0.65]{$3\times$} (0,1);
\draw[red, decoration={markings, mark=at position 0.5 with {\arrow{>}}}, postaction={decorate}] (1,0) --node[midway, left, red, scale=0.65]{$3\times$} (1,1);
\draw[red, decoration={markings, mark=at position 0.5 with {\arrow{>}}}, postaction={decorate}] (0,0) --node[midway, above, red, scale=0.65]{$3\times$} (1,0);
\foreach \x in {-1,0,1,2}{
	\foreach \y in {-1,0,1,2}{
		\def\clr{black}
		\ifnum \x=0 \ifnum \y=0 \def\clr{red} 
		\else \fi \else \fi
		\node[fill=\clr, circle, minimum size=0.4em, inner sep=0pt, outer sep=0pt] at (\x,\y){};
	}
}
\end{scope}	
\end{tikzpicture}
\caption{A visual representation of $S$ and $BS$.}
\label{F:square_ev}
\end{figure}

Accounting for changing orientations with signs, this means we can write the corresponding operator $B=TPI$ with respect to $\cal B_1$ as the $\C[\Z^2]$-valued matrix
$$ B = \frac{1}{6}\begin{pmatrix}
x+x^{-1} & 1-x-y^{-1}+xy^{-1} \\
1-x^{-1}-y+x^{-1}y & y+y^{-1}
\end{pmatrix}.
$$
We can readily verify that for $q\in [0,1]$ indeed
\begin{align*}
B_q &= q\Id + (1-q)B = q \Id + \frac{1-q}{6} \left(2\Id-\Delta\right) = C_{1,q}^{-1} \Id - C_{1,q}^{-1}C_{2,q} \Delta,
\end{align*}
and thus $C_{1,q} B_q = \Id - C_{2,q} \Delta$.
Looking at the boundary of $\beta_\circlearrowright$ given by 
$$S = (1-x)\alpha_\uparrow + (y-1)\alpha_\rightarrow,$$  
it is an eigenstate of $B$ with eigenvalue $\frac{1}{6}\left(x+x^{-1}+y+y^{-1}-2\right)$, compare Figure~\ref{F:square_ev}.

Here, $x+x^{-1}+y+y^{-1} = 4P^{\Z^2}$, where $P^{\Z^2}$ can formally also be interpreted as the propagation operator of the uniform nearest neighbour random walk on the grid $\Cayley(\Z^2)$ (or, in this case, rather the $2$-cells of $\R^2$ with this chosen CW structure). 
We denote $\lambda = \frac{1}{6}(4P^{\Z^2}-2)$. 
A straight-forward computation shows that $S$ is an eigenstate to $B_q$ with eigenvalue
$$ \lambda_q = C_{1,q}^{-1}C_{2,q}\left(4P^{\Z^2}+[C_{2,q}^{-1}-4]\right). $$
Since $C_{2,q}(4+[C_{2,q}^{-1}-4]) = 1$, 
we can set $q' = 1-4C_{2,q} = \frac{4q-1}{2q+1}$ and can formally interpret  
$$  4C_{2,q}P^{\Z^2}+[1-4C_{2,q}] = P^{\Z^2}_{q'} $$
as the propagation operator of the corresponding $q'$-lazy random walk on $\Cayley(\Z^2)$. 
Note that for $q\in [\nicefrac{1}{4},1]$ we have $4C_{2,q}\in [0,1]$ and $q'\in [0,1]$ so this makes sense. 
In particular, 
$$\lambda_q = C_{1,q}^{-1} P^{\Z^2}_{q'} \quad \text{ and } \quad B_q S = C_{1,q}^{-1} P^{\Z^2}_{q'}S. $$
The return quantity $p_q(n)$ that we are interested in is given by $p_q(n) = p_{q,\alpha_\uparrow}(n)+p_{q,\alpha_\rightarrow}(n)$, where $p_{q,\alpha_\uparrow}(n)=\langle B_q^n\alpha_\uparrow, \alpha_\uparrow\rangle$ is the coefficient of $1\alpha_{\uparrow}$ in $B_q^n\alpha_\uparrow$ and similarly for $\alpha_\rightarrow$. 
By symmetry, $p_{q,\alpha_\uparrow}(n)=p_{q,\alpha_\rightarrow}(n)$ so that 
$$ p_q(n) = 2p_{q,\alpha_\uparrow}(n) = 2\langle B_q^n \alpha_\uparrow, \alpha_\uparrow\rangle. $$
We note from Figure~\ref{F:R2Z2CW_Moves} that 
\begin{align}
\label{Eq_BqWithC1q}
B_q(\alpha_\uparrow) 
&= q\alpha_\uparrow + \frac{1-q}{6}\left(x^{-1}S - S +2\alpha_\uparrow \right)
= C_{1,q}^{-1} \alpha_\uparrow + C_{1,q}^{-1}C_{2,q} (x^{-1}-1) S.
\end{align} 
Since the random walk is $\Z^2$-invariant, this yields
\begin{align*}
B_q^n(\alpha_\uparrow)
&=  C_{1,q}^{-1}B_q^{n-1}(\alpha_\uparrow) + C_{1,q}^{-1}C_{2,q} (x^{-1}-1)B_q^{n-1} S 
\end{align*}
Resolving this recursive formula we obtain
\begin{align*}
B_q^n(\alpha_\uparrow)
&=  C_{1,q}^{-n}\alpha_\uparrow + \sum_{k=0}^{n-1} C_{1,q}^{-n+k} C_{2,q} (x^{-1}-1)B_q^{k} S\\
&=  C_{1,q}^{-n} \left(\alpha_\uparrow + \sum_{k=0}^{n-1}  C_{2,q}  (x^{-1}-1)(P^{\Z^2}_{q'})^{k} S\right)
\end{align*}
In order to find the coefficient of $1\alpha_\uparrow$, we notice that 
\begin{align*}
\langle S, 1\alpha_\uparrow \rangle &= 1, \qquad
\langle x^{-1}S, 
 1\alpha_\uparrow \rangle = -1 \quad \text{ and } \quad
\langle gS, 1\alpha_\uparrow \rangle = 0\quad \text{ for } g\notin \{1,x^{-1}\}
\end{align*}
and therefore it follows that 
\begin{align*}
\langle (P^{\Z^2}_{q'})^{k} S, 1\alpha_\uparrow\rangle &= \langle (P^{\Z^2}_{q'})^{k}, 1-x^{-1}\rangle,
\\
\langle x^{-1}(P^{\Z^2}_{q'})^{k} S, 1\alpha_\uparrow\rangle &= \langle (P^{\Z^2}_{q'})^{k}, x-1\rangle.
\end{align*}
Using this we obtain
\begin{align*}
\frac{1}{2}p_q(n) 
&= \langle B_q^n \alpha_\uparrow, \alpha_\uparrow\rangle
\\
&= C_{1,q}^{-n} \left(1 + \sum_{k=0}^{n-1}  C_{2,q} \left\langle(P^{\Z^2}_{q'})^{k}, (x-1)-(1-x^{-1})\right\rangle \right) \\
&= C_{1,q}^{-n} \left(1 + \sum_{k=0}^{n-1}  C_{2,q} \left\langle(P^{\Z^2}_{q'})^{k}, x+x^{-1}-2\right\rangle \right).
\end{align*}
By symmetry, the coefficients of $(P^{\Z^2}_{q'})^{k}$ for $x$ and $x^{-1}$ agree, hence
\begin{align*}
\frac{1}{2}C_{1,q}^{n}p_q(n) 
&= 1 - 2C_{2,q} \sum_{k=0}^{n-1}   \left\langle(P^{\Z^2}_{q'})^{k}, 1-x\right\rangle 
\\
&= 1 - 2C_{2,q} \sum_{k=0}^{n-1}  \left( p^{\Z^2}_{q'}(k) - p^{\Z^2}_{q'}(e\xrightarrow{k} x) \right)
\end{align*}
where $p^{\Z^2}_{q'}(k)$ is the return probability of the $q'$-lazy nearest neighbour random walk on $\Z^2$ after $k$ steps and $p^{\Z^2}_{q'}(e\xrightarrow{k} x)$ the probability of the random walk to be at the vertex $x$ after $k$ steps.
If we write 
$\mathbb E^g_q(n)$ for the expected number of visits of the vertex $g$ in the first $n$ steps for the $q$-lazy nearest neighbour random walk on $\Z^2$ (counting the starting position for $\mathbb E^e_q(n)$, if $q=0$ we suppress it in notation), we can write this as
$$ \frac{1}{2}C_{1,q}^n p_q(n) = 1- 2C_{2,q}(\mathbb E^e_{q'}(n-1)-\mathbb E^x_{q'}(n-1)). $$
Notice that $q' = 1-4C_{2,q}$ implies that $2C_{2,q} = \frac{1-q'}{2}$. Hence,
$$ \frac{1}{2}C_{1,q}^n p_q(n) = 1-\frac{1-q'}{2}(\mathbb E^e_{q'}(n-1)-\mathbb E^x_{q'}(n-1)). $$

For $q<1$ large enough, we expect that $$ \mathbb E^e_{q'}(n-1)-\mathbb E^x_{q'}(n-1) \sim 1 - \Theta(n^{-1}) \qquad \text{for } n\to\infty.  $$
Plugging this back into the equation above, this would imply that 
$$  p_q(n)  \sim C^{-n}_{1,q}\left(1 + \Theta(n^{-1})\right) \qquad \text{for } n\to\infty. $$
Here, we can read off $b^{(2)}(d_2^*) = 1$, corresponding to the kernel of $d_2^*$ of $\cal NG$-dimension one, and the Novikov-Shubin invariant $\alpha_1(\R^2) = \alpha(d_2^*) = 2$. 
	\printbibliography 
\end{document}